\documentclass[12pt,oneside]{amsproc}
\usepackage[T2A]{fontenc}
\usepackage[cp1251]{inputenc}
\usepackage[russian]{babel}
\usepackage{amssymb}
\usepackage{amsthm}

\newtheorem{thm}{Теорема}
\newtheorem*{rem}{Замечание}
\newtheorem*{deff}{Определение}

\newtheorem*{lemma}{Лемма}

\newcommand{\opt}{\operatorname{t}}

\newcommand{\ve}{\varepsilon}
\newcommand{\la}{\lambda}
\renewcommand{\sl}{\sum\limits}

\newcounter{raz}
\newcounter{dva}
\newcommand{\razs}{\addtocounter{raz}{1}\par\vspace{3 pt}\noindent$\mathrm{\arabic{raz}) }$ }
\newcommand{\dvas}{\addtocounter{dva}{1}\par\vspace{3 pt}\noindent$\mathrm{\arabic{dva}) }$ }

\title{Регулярные и вполне регулярные дифференциальные операторы}
\author{Е.~А.~Ширяев, А.~А.~Шкаликов}

\thanks{Работа выполнена при поддержке Российского фонда фундаментальных исследований, грант
№~04-01-00712 и программы <<Ведущие научные школы>>, грант №
НШ-5247.2006.1.}
\begin{document}
\vphantom{} \maketitle

1. {\bf Регулярные операторы}. В работе \cite{Birkhoff} Биркгоф
выделил класс обыкновенных дифференциальных операторов, для
которых получил оценки ядер Грина и доказал теорему о разложении
по собственным функциям. Операторы из этого класса он назвал
\emph{регулярными}. Напомним определение Биркгофа в упрощенной
форме.

Пусть $L$ --- оператор, порожденный дифференциальным выражением
\begin{equation}
\label{ly}
l(y)=(-i)^ny^{(n)}(x)+p_2(x)y^{(n-2)}+\cdots+p_n(x)y,
\end{equation}
и $n$ линейно независимыми краевыми условиями вида
\begin{equation}
\label{bc_0} U_j(y)=\sum\limits_{s=0}^{n-1}\left (
a_{j,\,s}y^{(s)}(0)+b_{j,\,s}y^{(s)}(1)\right )=0, \quad
j=1,\ldots,n.
\end{equation}
Предполагаем, что коэффициенты $p_j(x)$, $j=1,\ldots,n$,
--- суммируемые комплексные функции на отрезке $[0,1]$. Считаем,
что $L$ действует в пространстве $L_2(0,\,1)$  и определен
равенством $L(y)=l(y)$  на области
$$
D(L)=\{y\mid y^{(s)}\in AC[0,\,1],\;s=0,1,\ldots,n-1,\;l(y)\in
L_2, U_j(y)=0,\;j=1,\ldots,n\}.
$$
Здесь $AC[0,\,1]$ --- пространство абсолютно непрерывных функций.

Перепишем краевые условия, выделив старшие производные
\begin{equation}
\label{bc} U_{j}(y):=a_j y^{(k_j)}(0)+b_j y^{(k_j)}(1)+
\sum\limits_{s=0}^{k_j-1}\left (
a_{j,\,s}y^{(s)}(0)+b_{j,\,s}y^{(s)}(1)\right ) =0, \quad
j=1,\ldots,n.
\end{equation}

Число $k_j$ ($j=1,\ldots ,n$, $0\leqslant k_j\leqslant n-1$)
назовем \emph{порядком краевого условия}, а число
$\varkappa=k_1+\ldots+k_n$ --- \emph{суммарным порядком}. Заменим
краевые условия их линейными комбинациями, при которых суммарный
порядок наименьший. Считаем, что получившиеся краевые условия
имеют вид \eqref{bc}, причем нумерация такова, что $n-1 \geqslant
k_1 \geqslant k_2\geqslant \ldots \geqslant k_n \geqslant 0$.
Тогда старшие линейные формы
$$
U_j^0(y)=a_jy^{(k_j)}(0)+b_jy^{(k_j)}(1),\quad j=1,\ldots,n,
$$
также линейно независимы (иначе суммарный порядок можно понизить),
что влечет $k_j>k_{j+2}$, $a_jb_j\not =0$.

\begin{deff}
При четном $n=2m$ оператор $L$ назовем регулярным по Биркгофу,
если  не равен нулю  определитель
$$\theta=
\begin{vmatrix}
a_1 \varepsilon_1 ^{k_1}& \ldots &a_1 \varepsilon_{m-1} ^{k_1}&
a_1\varepsilon_{m} ^{k_1}& b_1 \varepsilon_{m+1} ^{k_1}&\ldots&
b_1 \varepsilon_{n} ^{k_1}
\\
\cdot&\cdot&\cdot&\cdot&\cdot&\cdot&\cdot
\\
a_n\varepsilon_1 ^{k_n}& \ldots &a_n \varepsilon_{m-1} ^{k_n}&
a_n\varepsilon_{m} ^{k_n}& b_n \varepsilon_{m+1} ^{k_n}&\ldots&
b_n \varepsilon_{n} ^{k_n}
\end{vmatrix},$$
где $\varepsilon_k=\exp(2\pi(k-1)i/n)$, $k=1,\ldots,n,$
--- корни $n$-й степени из $1$.
При нечетном $n=2m-1$ требуем, чтобы помимо этого определителя был
также отличен от нуля определитель, получающийся заменой в $m-$ом
столбце чисел $a_j$ на $b_j$.
\end{deff}
При $n=2m$ приведенное определение отличается от классического
\cite{Birkhoff}, \cite{Naimark}, где требуется отличие от нуля еще
одного числа, но несложно показать, что второе число может
обращаться в нуль только одновременно с приведенным здесь первым
числом. Коэффициенты $p_j(x)$ дифференциального выражения и
младшие члены в краевых условиях не участвуют в определении
регулярности; $L$ регулярен тогда и только тогда, когда регулярен
оператор $L_0(y)=(-i)^ny^{(n)}$, порожденный краевыми условиями
$U_j^0(y)=0$, $j=1,\ldots,n$.

Определение регулярности переносится на существенно более общие
классы краевых задач со спектральным параметром. Для обыкновенных
дифференциальных уравнений и систем это сделано в \cite{Tamarkin},
\cite{Sh1}, \cite{Lu}, а теория регулярных  задач для уравнений с
частными производными построена в работах \cite{Lo}, \cite{AN},
\cite{AV}. Интересно отметить, что для обыкновенных
дифференциальных операторов порядка $2m$  в случае распадающихся
краевых условий можно определить регулярность по Лопатинскому и
оно будет эквивалентным регулярности по Биркгофу
(см.~\cite[\S~9]{Sh2}).

Для формулировки основной теоремы о регулярных операторах удобно
ввести дополнительные обозначения. Лучи $\arg\rho=\pm(\pi/2n
+\arg(i\varepsilon_k))$ в комплексной $\rho-$плоскости назовем
критическими. Их число равно $n$ и $2n$ в четном и нечетном
случаях, соответственно. Вырежем из $\rho$-плоскости открытые
секторы раствора $\varepsilon<\pi/(2n)$, биссектрисами которых
служат критические лучи, а оставшиеся $n$ или $2n$ замкнутых
секторов обозначим через $\Omega(\ve)$.

Если резольвентное множество оператора $L$  непусто, то из
компактности вложения $D(L)\subset L_2(0,\,1)$ (мы снабжаем $D(L)$
нормой графика оператора $L$, тогда $D(L)$  становится
гильбертовым пространством) следует, что спектр $L$ дискретный,
т.е. состоит из изолированных собственных значений
$\{\lambda_k\}_1^\infty$ конечной алгебраической кратности.
Обозначим через $\rho_{k,\,j}$, $k=1,\ldots,n$,  корни уравнения
$\la_j=\rho^n$, через $B_{k,\,j}(\delta)$ --- круги радиуса
$\delta$ в $\rho$-плоскости с центрами в точках $\rho_{k,\,j}$, и
через $B(\delta)$ --- объединение всех таких кругов по обоим
индексам $k$ и $j$. Известно \cite{Naimark}, что резольвента $L$
есть интегральный оператор
$$
(L-\rho^n)^{-1}f(x)=\int_0^1G(x,\xi,\rho)f(\xi)d\xi,
$$
где $G(x,\xi,\rho)$ --- ядро Грина. Теперь сформулируем основной
результат.
\begin{thm}
Следующие утверждения эквивалентны. \razs Оператор $L$ регулярен
по Биркгофу; \razs при любом $\delta>0$ всякий некритический луч,
выходящий из нуля, асимптотически не пересекает множество $\Bbb C
\backslash B(\delta)$, и для всех $x,\,\xi\in[0,1]$ и $\rho\in
\Bbb C \backslash B(\delta)$ для функции Грина справедлива  оценка
\begin{equation}\label{|G|}
|G(x,\,\xi,\,\rho)|\leqslant M |\rho|^{-n+1},
\end{equation}
 где постоянная $M=M(\delta)$ не зависит от $x, \xi, \rho$;
\razs найдется последовательность точек $\{\rho_k\}$ в одном из
секторов множества $\Omega(\ve)$ (в случае нечетного $n$
--- две последовательности в двух соседних секторах этого
множества), такая, что $\rho_k\to\infty$ и оценка \eqref{|G|}
выполнена при $\rho =\rho_k$; \razs при любом $\ve>0$ в области
$\Omega(\ve)$  выполняется оценка
 \begin{equation}\label{||G||}
 \|(L-\rho^n)^{-1}\|\leqslant M |\rho|^{-n},\quad |\rho|
 \geqslant |\rho_0|,
 \end{equation}
где постоянная $M=M(\ve,\rho_0)$ не зависит от $\rho$, а
$\|\cdot\|$ означает норму в $L_2$; \razs найдется
последовательность точек $\{\rho_k\}$ в одном из секторов
множества $\Omega(\ve)$ (в случае нечетного $n$
--- две последовательности в двух соседних секторах этого
множества), такая, что $\rho_k\to\infty$ и оценка \eqref{||G||}
выполнена при $\rho =\rho_k$; \razs система собственных и
присоединенных функций оператора $L$ образует безусловный базис со
скобками в пространстве $L_2$, причем в скобки объединяются не
более двух собственных функций.
\end{thm}

\begin{proof}  Импликация $1)\Rightarrow 2)$  доказана Биркгофом
\cite{Birkhoff}. Формально в \cite{Birkhoff} (см. также
\cite{Naimark}) оценка \eqref{|G|} доказана на некоторой
последовательности контуров, уходящих на бесконечность. Однако при
условии регулярности, оценку снизу характеристического
определителя (голоморфной функции, нули которой совпадают с
числами $\rho_{k,\,j}$) легко провести вне множества $B(\rho)$.
Тогда оценка \eqref{|G|} получается при всех $\rho\in\mathbb{C}
\backslash B(\rho)$. Импликация $1)\Rightarrow 4)$, доказана
Бензингером \cite{Benz}. Импликации $2)\Rightarrow 3)$ и
$4)\Rightarrow 5)$ тривиальны, а $5)\Rightarrow 1)$ доказана
Минкиным в его недавней работе \cite{Min}. Импликация
$1)\Rightarrow 6)$ доказана Шкаликовым \cite{Sh3} (там же см.
ссылки на предшествующие работы Н.~Данфорда и Дж.~Шварца,
Г.~М.~Кесельмана и В.~П.~Михайлова на эту тему). Остается доказать
$4)\Rightarrow 1)$ и $6)\Rightarrow 1)$. Доказательства обеих этих
импликаций нетривиальны, в них существенно используются результаты
и приемы недавней работы Минкина \cite{Min}. Подробное изложение
будет дано в другой работе, здесь отметим только один
вспомогательный результат, играющий важную роль.

\begin{lemma}
Пусть $\{\alpha_k\}_1^\infty$ --- комплексные, не равные нулю
числа, аргументы которых различны. Положим
$$
\psi_j(x)=\sum_{k=1}^nc_{k,\,j}e^{\alpha_k\rho_jx}(1+\varphi_{k,\,j}(x)),
$$
где $c_{k,\,j}$ --- произвольные числа, такие что
$|c_{1,\,j}|+\ldots+|c_{n,\,j}|\not =0$, а
$\|\varphi_{k,\,j}\|^2=o(|\rho_j|^{-1})$ при $|\rho_j|\to\infty$.
Назовем критическими лучи, аргументы которых равны $\pm
\pi/2+\arg\alpha_k$. Последовательность $\{\rho_j\}_1^\infty$,
занумерованную в порядке возрастания модулей, назовем редкой, если
найдется целое число $l$, такое, что $|\rho_{j+l}/\rho_j|\geqslant
2$ при всех $j\geqslant 1$.

Пусть система  $\{\psi_j(x)\}_1^\infty$ образует базис Рисса со
скобками в замыкании своей линейной оболочки в пространстве
$L_2(0,\,1)$, причем в скобки заключаются не более фиксированного
числа $N$  функций. Тогда в каждом замкнутом секторе, не
содержащем критические лучи, последовательность $\{\rho_j\}$
является редкой.
\end{lemma}
\end{proof}
%\vspace{0.5 cm}
\noindent2. {\bf Вполне регулярные операторы}. \noindent Здесь
будем рассматривать дифференциальные операторы четного порядка
$n=2m$, заданные дифференциальным выражением
\begin{equation}
\label{l_diver} l(y) = \sum_{k=0}^m (-1)^k \left\{\left(
p_k(x)y^{(k)}\right)^{(k)} - \left[ \left(
q_k(x)y^{(k)}\right)^{(k-1)}+\left(
r_k(x)y^{(k-1)}\right)^{(k)}\right]\right\},
\end{equation}
где $p_m(x)=1$, $r_0(x)=0$. Чтобы не осложнять существо дела,
предположим, что коэффициенты $p_k(x)$, $q_k(x)$, $r_k(x)$ таковы,
что после раскрытия производных дифференциальное выражение
приводится к виду \eqref{ly}. Для этого достаточно, чтобы
$p_k(x),\,r_k(x)\in W_1^k[0,1]$, $q_k(x)\in W_1^{k-1}[0,1]$.

Введем квазипроизводные
\begin{equation}
\label{qua}
\begin{array}{r}
y^{[k]} =y^{(k)}, \quad k=0,\dots ,m-1, \quad y^{[m]}=y^{(m)}-r_my^{(m-1)},\\
y^{[m+k]}=-(y^{[m+k-1]})'+p_{m-k}y^{(m-k)}+\left[
q_{m-k+1}y^{(m-k+1)}-r_{m-k}y^{(m-k-1)}\right],\\ k=1,\dots,m.
\end{array}
\end{equation}
Легко видеть, что при фиксированном $x$ функционал $y^{[k]}(x)$
является линейной комбинацией функционалов $y^{(s)}(x)$,
$s=0,\ldots,k$, и наоборот. Поэтому краевые условия можно записать
в виде
\begin{equation}
\label{bc_q}
By^\land+Cy^\lor=0,
\end{equation}
где $B$ и $C$ --- некоторые матрицы, а $y^\land$ и $y^\lor$
--- векторы значений квазипроизводных:
\begin{flushleft}
\quad\quad$y^{\land}=\bigl(y(0),\; y'(0),\ldots,
y^{(m-1)}(0),\;y(1),\;
y'(1),\ldots, y^{(m-1)}(1) \bigr)^{\opt},$\\
\quad\quad$y^{\lor}=\bigl(y^{[2m-1]}(0),\; y^{[2m-2]}(0),\ldots,
y^{[m]}(0),\;-y^{[2m-1]}(1),\; -y^{[2m-2]}(1),\ldots, -y^{[m]}(1)
\bigr)^{\opt},$
\end{flushleft}
где верхний индекс $\opt$ означает транспонирование.
\begin{deff}
Оператор $L$, порожденный дифференциальным выражением
\eqref{l_diver} и краевыми условиями \eqref{bc_q}, назовем вполне
регулярным, если выполнено условие
${B^{-1}(\operatorname{im}C)=\mathbb{C}^{2m}\ominus\operatorname{ker}{C}}$,
где $B^{-1}$ понимается как взятие полного прообраза от
$\operatorname{im}C$ при отображении $B$.
\end{deff}
\begin{thm}
Пусть оператор $L$ задан выражением \eqref{l_diver} и краевыми
условиями \eqref{bc_q}. Следующие утверждения эквивалентны . \dvas
Оператор $L$ является вполне регулярным. \dvas Числовой образ
оператора $L$ не совпадает со всей комплексной плоскостью. \dvas
Числовой образ оператора $L$ содержится в некоторой полуплоскости
в $\mathbb{C}$. \dvas Квадратичная форма оператора
$(Ly,\,y)_{L_2}$ представима в виде
\begin{equation}
\label{Lyy}
\sl_{k=1}^m\left[(p_ky^{(k)},\,y^{(k)})_{L_2}+(q_ky^{(k)},\,y^{(k-1)})_{L_2}-(r_ky^{(k-1)},\,y^{(k)})_{L_2}\right]
+(Ay^\land,\,y^\land)_{\mathbb{C}^{2m}}.
\end{equation}
\end{thm}
\begin{proof}
Эквивалентность утверждений 2) и 3) следует из свойства выпуклости
числового образа.

Импликация 1)$\Rightarrow$4) для случая $q_k=r_k\equiv 0$ доказана
в работе \cite{AAB}. В общем случае доказательство можно провести
аналогично. Чтобы показать справедливость обратного утверждения
докажем следующую лемму.
\begin{lemma}
Пусть матрица $A$ такова, что для всех $x$, $y_1$,
$y_2\in\mathbb{C}^{2m}$ справедливо равенство
$(y_2,x)_{\mathbb{C}^{2m}}=(Ay_1,x)_{\mathbb{C}^{2m}}$, если $x\in
B^{-1}(\operatorname{im}C)$ и $By_1+Cy_2=0$. Тогда оператор $L$
вполне регулярен.
\end{lemma}
\begin{proof}
Сначала покажем, что
$B^{-1}(\operatorname{im}C)\perp\operatorname{ker}C$. Пусть
векторы $y_1$ и $y_2$ таковы что $By_1+Cy_2=0$. Тогда $y_1\in
B^{-1}(\operatorname{im}C)$ и
${(y_2,\,y_1)_{\mathbb{C}^{2m}}}=(Ay_1,\,y_1)_{\mathbb{C}^{2m}}$.
Поскольку ${By_1+C(y_2+v)=0}$ для любого
$v\in\operatorname{ker}C$, то
$(y_2+v,\,y_1)_{\mathbb{C}^{2m}}=(Ay_1,\,y_1)_{\mathbb{C}^{2m}}$.
Поэтому  $(v,\,y_1)_{\mathbb{C}^{2m}}=0$, т.е.
$B^{-1}(\operatorname{im}C)\perp\operatorname{ker}C$.

Равенство
${B^{-1}(\operatorname{im}C)=\mathbb{C}^{2m}\ominus\operatorname{ker}C}$
следует теперь из того, что краевые условия линейно независимы.
\end{proof}
Доказательство импликации 4)$\Rightarrow$3) следует из
компактности вложений $W_2^{r}\subset W_2^{m}$ при  $r\leqslant
m-1$ и компактности операторов следа $T_0y=y^{(r)}(0)$ и
$T_1y=y^{(r)}(1)$ как операторов из $W_2^m$ в $\mathbb{C}$.
Доказательство обратной импликации более сложно и здесь мы его не
приводим. В нем используется прием из работы \cite{Eshir}.
\end{proof}
\begin{rem}
Всякий вполне регулярный оператор регулярен. Действительно, если
$L$  вполне регулярен, то его числовой образ лежит в некоторой
полуплоскости, а потому найдется сектор в комплексной плоскости, в
котором резольвента $(L-\rho^m)^{-1}$ имеет оценку \eqref{||G||}.
Тогда по теореме 1 оператор $L$ регулярен. Обратное, вообще
говоря, неверно. Это показывает следующий пример.
\end{rem}
\noindent{\bfseries Пример.} Пусть $Ly=y^{(4)}$, а краевые условия
имеют вид
$$
\left \{
\begin{array}{r}
-y'''(0)+y''(0)+y(0)=0\\
y'''(1)+y(1)=0\\
y'(0)=0\\
y'(1)=0
\end{array}
\right ..
$$
Простая проверка показывает, что $L$ регулярен, но не вполне
регулярен.

\end{document}